\documentclass[a4paper,12pt]{amsart}
\usepackage{amssymb}
\usepackage{ifthen}
\usepackage{graphicx}
\usepackage{color,xcolor}
\nonstopmode
\numberwithin{equation}{section}
\newtheorem{thm}{Theorem}
\newtheorem{cor}{Corollary}
\newtheorem{lem}{Lemma}
\newtheorem{prop}{Proposition}

\newtheorem{conj}{Conjecture}
\newtheorem{prob}{Problem}

\theoremstyle{definition}
\newtheorem{defn}{Definition}
\newtheorem{ca}{Case}

\newtheorem{rem}{Remark}

\newenvironment{pf}[1][]{%
 \vskip 1mm
 \noindent
 \ifthenelse{\equal{#1}{}}%
  {{\slshape Proof. }}%
  {{\slshape #1.} }%
 }%
{\qed\medskip}

\newcounter{alphabet}

\newenvironment{Thm}[1][]{\refstepcounter{alphabet}%
\bigskip%
\noindent%
{\bf Theorem \Alph{alphabet}}%
\ifthenelse{\equal{#1}{}}{}{ (#1)}%
{\bf .} \itshape}{\vskip 8pt}

\newenvironment{Lem}[1][]{\refstepcounter{alphabet}%
\bigskip%
\noindent%
{\bf Lemma \Alph{alphabet}}%
{\bf .} \itshape}{\vskip 8pt}

\newcounter{alphabet2}

\newcommand{\IN}{{\mathbb N}}
\newcommand{\IC}{{\mathbb C}}
\newcommand{\ID}{{\mathbb D}}

\def\be{\begin{equation}}
\def\ee{\end{equation}}

\newcommand{\ben}{\begin{enumerate}}
\newcommand{\een}{\end{enumerate}}

\newcommand{\blem}{\begin{lem}}
\newcommand{\elem}{\end{lem}}
\newcommand{\bthm}{\begin{thm}}
\newcommand{\ethm}{\end{thm}}
\newcommand{\bcor}{\begin{cor}}
\newcommand{\ecor}{\end{cor}}
\newcommand{\beg}{\begin{exam}}
\newcommand{\eeg}{\end{exam}}
\newcommand{\begs}{\begin{examples}}
\newcommand{\eegs}{\end{examples}}
\newcommand{\bdefe}{\begin{defn}}
\newcommand{\edefe}{\end{defn}}
\newcommand{\bprob}{\begin{prob}}
\newcommand{\eprob}{\end{prob}}
\newcommand{\bques}{\begin{ques}}
\newcommand{\eques}{\end{ques}}
\newcommand{\bei}{\begin{itemize}}
\newcommand{\eei}{\end{itemize}}
\newcommand{\bcon}{\begin{conj}}
\newcommand{\econ}{\end{conj}}
\newcommand{\bop}{\begin{op}}
\newcommand{\eop}{\end{op}}

\newcommand{\bas}{\begin{assertion}}
\newcommand{\eas}{\end{assertion}}

\newcommand{\bfa}{\begin{fact}}
\newcommand{\efa}{\end{fact}}

\newcommand{\bca}{\begin{ca}}
\newcommand{\eca}{\end{ca}}

\newcommand{\bst}{\begin{step}}
\newcommand{\est}{\end{step}}

\newcommand{\bsca}{\begin{sca}}
\newcommand{\esca}{\end{sca}}

\newcommand{\bcl}{\begin{cl}}
\newcommand{\ecl}{\end{cl}}

\newcommand{\bmlem}{\begin{mlem}}
\newcommand{\emlem}{\end{mlem}}

\newcommand{\bscl}{\begin{scl}}
\newcommand{\escl}{\end{scl}}

\newcommand{\bcons}{\begin{conjs}}
\newcommand{\econs}{\end{conjs}}

\newcommand{\bprop}{\begin{prop}}
\newcommand{\eprop}{\end{prop}}

\newcommand{\br}{\begin{rem}}
\newcommand{\er}{\end{rem}}
\newcommand{\brs}{\begin{rems}}
\newcommand{\ers}{\end{rems}}
\newcommand{\bo}{\begin{obser}}
\newcommand{\eo}{\end{obser}}
\newcommand{\bos}{\begin{obsers}}
\newcommand{\eos}{\end{obsers}}
\newcommand{\bpf}{\begin{pf}}
\newcommand{\epf}{\end{pf}}
\newcommand{\ba}{\begin{array}}
\newcommand{\ea}{\end{array}}
\newcommand{\beq}{\begin{eqnarray}}
\newcommand{\beqq}{\begin{eqnarray*}}
\newcommand{\eeq}{\end{eqnarray}}
\newcommand{\eeqq}{\end{eqnarray*}}

\newcommand{\ra}{\to}

\newcommand{\ds}{\displaystyle}

\newcounter{minutes}
\divide\time by 60
\newcounter{hours}
\multiply\time by 60 \addtocounter{minutes}{-\time}

\begin{document}

\bibliographystyle{amsplain}
\title [Multidimensional analogues of refined Bohr's inequality]
{Multidimensional analogues of refined Bohr's inequality}

\def\thefootnote{}
\footnotetext{ \texttt{\tiny File:~\jobname .tex,
          printed: \number\day-\number\month-\number\year,
          \thehours.\ifnum\theminutes<10{0}\fi\theminutes}
} \makeatletter\def\thefootnote{\@arabic\c@footnote}\makeatother


\author{Ming-Sheng Liu
}
 \address{M-S Liu, School of Mathematical Sciences, South China Normal University, Guangzhou, Guangdong 510631, China.} \email{liumsh65@163.com}

\author{Saminathan Ponnusamy}
\address{S. Ponnusamy, Department of Mathematics,
Indian Institute of Technology Madras, Chennai-600 036, India. }
\email{samy@iitm.ac.in}




\begin{abstract}
In this paper, we first establish a version of multidimensional analogues of the refined Bohr's inequality.
Then we establish two versions of multidimensional analogues of improved Bohr's inequality with initial coefficient
being zero. Finally we establish two versions of multidimensional analogues of improved Bohr's inequality with the initial
coefficient being replaced by absolute value of the function,  and
to prove that most of the results are sharp. 
\end{abstract}

\maketitle
\pagestyle{myheadings}
\markboth{M-S Liu and S. Ponnusamy}{Multidimensional analogues of refined Bohr's inequality}

\section{Preliminaries and some basic questions}\label{HLP-sec1}

The classical theorem of Bohr \cite{B1914}, examined a century ago,  generates intensive research activities on Bohr's phenomena
and it was revived by many with great interest in the nineties due to the extensions to holomorphic functions of several complex variables and to more abstract settings.
For example in 1997, Boas and Khavinson \cite{BK1997} defined $n$-dimensional Bohr radius for the family of holomorphic functions bounded by $1$ on the unit polydisk
(see Section \ref{sect1-3} for details). This paper stimulated interests on Bohr type questions in different settings.
For example, Aizenberg \cite{A2000,A2005,A2007}, Aizenberg et al. \cite{AAD2000,AAD,AV2004,AV2008}, Defant and Frerick \cite{DF}, and
Djakov and Ramanujan \cite{DjaRaman-2000} have established further results on Bohr's phenomena for multidimensional
power series. Several other aspects and generalizations of Bohr's inequality may be obtained from
\cite{AAP2016, AlKayPON-19,DFOOS,HHK2009,KSS2017,LP2019,PPS2002,PS2004,PS2006},
the monograph of Kresin and  Maz'ya \cite{KM2007}, and the references therein. In particular,  \cite[Section 6.4]{KM2007} on Bohr's type theorems
contains rich opportunities to extend several inequalities to holomorphic functions of several complex variables and more importantly to
solutions of partial differential equations.

 In particular, after the appearance of some recent articles  \cite{AliBarSoly,KP2017,KP2018, KayPon3},
several new problems on Bohr's inequality in the plane case are investigated (cf. \cite{BhowDas-18,BhowDas-19,KayPon3, LSX2018, LPW2020,PVW2019,PVW201911}).


Our primary interest in this paper is to establish several multidimensional analogues of improved Bohr's inequality for the analytic case
and to prove that most of the results are sharp.

\subsection{Classical Inequality of H. Bohr}
Let us fix some notations. Throughout the discussion, let $H_\infty $ denote the class of all bounded analytic functions $f$ on the unit disk $\ID$
with the supremum norm $\|f\|_\infty :=\sup_{z\in \ID}|f(z)|$,
$${\mathcal B} = \{f\in H_\infty :\, \|f\|_\infty \leq 1 \} ~\mbox{ and }~
{\mathcal B}_0=\{\omega \in {\mathcal B}:\, \omega (0)=  0 \}.
$$
Then the classical inequality examined by Bohr in 1914 \cite{B1914} states that $1/3$ is the largest value of $r\in [0,1)$ for which
the following inequality holds:
\begin{equation}
B(f,r):=\sum_{k=0}^{\infty}|a_k|r^k\leq 1
\label{liu1}
\end{equation}
for every analytic function $f\in \mathcal{B}$ with the Taylor series expansion $f(z)=\sum_{k=0}^{\infty}a_kz^k$.
Bohr actually obtained that  \eqref{liu1} is true when  $r\leq 1/{6}$. Later  M.~Riesz, I.~Schur and  F. Wiener, independently
established the Bohr inequality \eqref{liu1} for $r\leq 1/{3}$ and that $1/3$
is the best possible constant. It is quite natural that the constant $1/3$ is called the Bohr radius
for the space $\mathcal{B}$. 
Moreover, for
$$\varphi_a(z)=\frac{a-z}{1-a z},\quad a\in [0,1),
$$
it follows easily that $B(\varphi_a,r)>1$ if and only if $r>1/(1+2a)$, which for $a\ra 1$ shows that $1/3$ is optimal.
It is worth pointing out that there is no extremal function in $\mathcal{B}$ such that the Bohr radius is precisely $1/3$ (cf. \cite[Corollary 8.26]{GarMasRoss-2018}).

Later it was shown by  Paulsen et al.   \cite[Corollary 2.7]{PPS2002} that if the constant term  $|a_0|$ in \eqref{liu1} is replaced by $|a_0|^2$, then the Bohr radius $1/3$ could be replaced by $1/2$ which is also optimal. Later in 1962, Bombieri \cite{Bom-62} (see also \cite{PPS2002}) found that the Bohr radius for the family ${\mathcal B}_0$ is $1/\sqrt{2}$ which is again optimal as the  function $z\varphi_a(z)$ demonstrates, where $\varphi_a(z)$ is defined as above.
See \cite{KP2018,KayPon3,KayPon_AAA18} for new proofs of it in a general form.

In \cite{BDK5},  B\'en\'eteau et al. proved that there is no Bohr phenomenon in Hardy spaces $H_p$ ($0<p<\infty)$ with the usual norm. On the other hand,
$p$-Bohr phenomenon for series of the type $\sum_{k=0}^{\infty}|a_k|^pr^k$ ($1\leq p$), for  $f(z)=\sum_{k=0}^{\infty}a_kz^k \in H_\infty $ was investigated first by Djakov and Ramanujan \cite{DjaRaman-2000} and in \cite{BDK5}. More recently, this was discussed again  by Kayumov and Ponnusamy \cite{KayPon_AAA18}, and the corresponding conjecture about Bohr radius was settled. 
In \cite{DjaRaman-2000}, the authors extended the notion of $p$-Bohr radius to the case of holomorphic functions of
several variables. In \cite{BDK5}, the authors considered Bohr type inequality in fairly general normed spaces of analytic functions and extended certain results to
the case of several variables.

\subsection{The refined versions of classical Bohr's inequality}

For $f(z)=\sum_{n=0}^{\infty} a_{n} z^{n}\in {\mathcal B}$, we let for convenience $f_0(z):=f(z)-f(0)$,
$$\|f_0\|_r^2  :=   \sum_{n=1}^{\infty}\left|a_{n}\right|^{2} r^{2 n} ~\mbox{ and }~
 B_{N}(f,r) :=  \sum_{n=N}^{\infty} |a_n| r^n  ~\mbox{ for $N\geq 0$},
$$
so that $B_0(f,r) = |a_0|+ B_1(f,r)$.

Recently, Kayumov and Ponnusamy \cite{KP2017,KP2018}, and Ponnusamy et al. \cite{PVW2019,PVW201911} established several refined
versions and improved versions of Bohr's inequality in the case of bounded analytic functions. We now recall a couple of them.

\begin{Thm}\label{Theo-A}
Suppose that $f \in\mathcal{B}$ has the expansion $f(z)=\sum_{n=0}^{\infty} a_{n} z^{n}$,  $f_0(z)=f(z)-f(0)$,
and $p>0$. Then 
we have the following:
\begin{enumerate}
\item[{\rm (a)}] 
$\ds B_1(f,r)+\left(\frac{1}{1+|a_0|}+\frac{r}{1-r}\right)\|f_0\|_r^2 \leq \frac{r}{1-r}(1-|a_0|^2)
$
\item[{\rm (b)}] 
$\ds
|a_0|^p+B_1(f,r)+ \left(\frac{1}{1+|a_0|}+\frac{r}{1-r}\right)\|f_0\|_r^2 \leq 1 ~\mbox{ for $\ds r\leq r_p(|a_0|)$},
$
where
$$r_p(|a_0|)=\frac{1-|a_0|^p}{2-|a_0|^2-|a_0|^p}.
$$
The numbers $r_p(|a_0|)$ and $\frac{1}{1+|a_0|}$ cannot be improved.
In particular, $r_1(|a_0|)=\frac{1}{2+|a_0|}$ and $r_2(|a_0|)=1/2$. Moreover,
$$\inf _{|a_0|<1} r_p(|a_0|) =\frac{p}{2+p} ~\mbox{ for $0< p\leq 2$}.
$$

\item[{\rm (c)}] 
If $f(0)=0$, then $B_1(f,r) \leq 1$ for $r\leq \frac{1}{\sqrt{2}}$, and the number $\frac{1}{\sqrt{2}}$ is sharp.
\end{enumerate}
\end{Thm}
\bpf
The proof of the first part (a) is obtained from \cite{PVW201911} while the proof of the second part (b) is remarked in \cite[Remark 1]{PVW2019}.
The proof of the third part (c) is due to Bombieri \cite{Bom-62}. See \cite{KP2017,KP2018} where
one can find a more general result for $p$-symmetric bounded functions. Further related results on this topic
may be found from \cite{DjaRaman-2000} and \cite[Corollary 2.9]{PPS2002}.
%
%
%
\epf

%

\begin{Thm} \label{Theo-C}(\cite{KP201701}) \, 
Suppose that $f\in \mathcal{B}$ and $f(z)=\sum_{n=0}^{\infty}a_{n} z^{n}$. Then for $p=1,2$,
\begin{equation}
|f(z)|^p+B_N(f,r) \leq 1~ \mbox{ for }~ r\leq R_{N,p},
\label{Liu15-1}
\end{equation}
where $R_{N,p}$ is the positive root of the equation $2(1+r)r^{N}-p(1-r)^{2}=0$. The radius $R_{N,p}$ is the best possible.
\end{Thm}

In fact in Section \ref{HLP-sec3}, we shall show that the sharp inequality \eqref{Liu15-1} actually continues to hold for all $p\in (0,2]$.

\begin{Thm}\label{Theo-D}(\cite{LLP2020}) \,
Suppose that $f\in \mathcal{B}$ and $f(z)=\sum_{n=0}^{\infty}a_{n} z^{n}$. Then for $p=1,2$,
\begin{equation}
|f(z)|^p+B_1(f,r)+\left(\frac{1}{1+|a_0|}+\frac{r}{1-r}\right)\|f_0\|_r^2\leq 1\quad for\quad |z|=r\leq r_{|a_0|, p},
\label{Liu015}
\end{equation}
where $r_{a, p}$ is the minimum positive root in $(0, 1)$ of the equation $A_{a,p}(r)=0$, where $a=|a_0|$, and
\begin{equation}%
A_{a,p}(r) =[1-(2-a^2)r](1+a r)^p - (1-r)(r+a)^p.
\label{Liu34}
\end{equation}
The radius $r_{|a_0|, p}$ is best possible. 
\end{Thm}

In fact in Section \ref{HLP-sec3}, we shall show that the sharp inequality \eqref{Liu015} actually continues to hold for all $p>0$.


\subsection{Multidimensional analogues of Bohr's inequality}\label{sect1-3}


Let $\mathbb{D}^n=\{z\in \IC^n: \, z=(z_1,\ldots, z_n), ~|z_j|<1,\, j=1,\ldots, n\}$ be the open unit polydisk.
Denote by  $\IN _0^n$ the set of all  $n$-tuple of non-negative integers. Also, denote by $K_n$ the largest non-negative number
such that if the $n$-variable power series
\begin{eqnarray}
\sum\limits_{\alpha} c_{\alpha}\, z^{\alpha}
\label{Liu17}
\end{eqnarray}
converges in   $\mathbb{D}^n$ and its sum $f$ has modulus less than $1$ so that $c_{\alpha}=\partial ^\alpha f(0)/\alpha !$, then
\begin{eqnarray}
\sum\limits_{\alpha} |c_{\alpha}|\, |z^{\alpha}|<1 ~\mbox{ for all $z\in K_n\cdot\mathbb{D}^n$.}
\label{Liu19}
\end{eqnarray}
Here $z=(z_1,\ldots, z_n)\in \IC^n$ is an $n$-tuple of complex numbers,
$$|\alpha|=|\alpha_1| +\cdots +|\alpha_n|,~~ \alpha !=\alpha_1! \cdots  \alpha_n! ~\mbox{ and }~
z^\alpha=z_1^{\alpha_1}\cdots z_n^{\alpha_n}.
$$

In 1997, Boas, Khavinson \cite{BK1997} obtained the following result as a multidimensional generalization of Bohr's inequality.

\begin{Thm}\label{Theo-DE}(\cite{BK1997})
It is true for $n>1$ that 
\begin{equation}
\frac{1}{3\sqrt{n}}<K_n<\frac{2\sqrt{\log n}}{\sqrt{n}}.
\label{Liu110}
\end{equation}
\end{Thm}

Clearly, Theorem~D 
shows that the Bohr radius $K_n$ for the case of polydisk depends on the
dimension $n$, and $K_n$ tend to zero as the dimension $n$ of the domain increases to $+\infty$.
For the unit hypercone $D^{o}= \{z=(z_1,\ldots, z_n)\in \IC^n:\,\sum_{j=1}^n |z_j|<1 \}$, Aizenberg \cite{A2000} obtained estimates that do not depend on $n$. In
\cite{A2000}, one can also find results for Reinhardt domain and complete bounded circular domain.
In \cite{AAD2000}, existence of Bohr's phenomenon under very general conditions is established.

\bdefe
A domain $\Omega\subset\mathbb{C}^n$ is called a Reinhardt domain centered at $0\in \Omega$ if for any
$z=(z_1,\ldots,z_n)\in \Omega$, and for each $\theta_k\in [0, 2\pi],\, k=1,\ldots, n$,
we have that $(z_1e^{i\theta_1},\ldots,z_n e^{i\theta_n})\in \Omega$. Furthermore, if for each $z=(z_1,\ldots,z_n)\in \Omega$,
and for each $|\xi_k|\leq 1,\, k=1,\ldots, n$, we have that $(z_1\xi_1,\cdots,z_n \xi_n)\in \Omega$, then $\Omega$ is called a complete Reinhardt domain.

A domain $Q\subset\mathbb{C}^n$ is a circular domain centered at $0\in Q$ if for any $z=(z_1,\ldots,z_n)\in Q$, and for each $\theta\in [0, 2\pi]$, we have that $e^{i\theta} z=(z_1e^{i\theta},\ldots,z_n e^{i\theta})\in Q$. A domain $Q\subset\mathbb{C}^n$ is called a complete circular domain centered at $0\in Q$ if for any $z\in Q$, and for each $\xi\in \overline{\mathbb{D}}$, we have that $ \xi\cdot z = (z_1\xi,\ldots,z_n \xi)\in Q$. For example, the balls and polydiscs in $\mathbb{C}^n$ are complete Reinhardt domains and complete circular domains.
\edefe

If $Q$ is a complete circular domain centered at $0\in Q\subset\mathbb{C}^n$, then every holomorphic function $f$ in $Q$
can be expanded into homogeneous polynomials given by
\begin{equation} \label{Liu13}
f(z)= \sum_{k=0}^{\infty}P_k(z) ~\mbox{ for $z\in Q$},
\end{equation}
where $P_k(z)$ is a homogeneous polynomial of degree $k$, and $P_0(z)=f(0)$.

In 2000, Aizenberg \cite{A2000} obtained the following multidimensional analogues of Bohr's inequality.

\begin{Thm}\label{Theo-E}(\cite{A2000})
If the series \eqref{Liu13} converges in the domain $Q$ and the estimate $|f(z)|<1$ holds in it, then
\begin{equation}
\sum_{k=0}^{\infty}|P_k(z)|<1
\label{Liu14}
\end{equation}
in the homothetic domain $(1/3)Q$. Moreover, if $Q$ is convex, then $1/3$ is the best possible constant.

\end{Thm}


It is natural to raise the following.

\bprob\label{HLP-prob1}
Can we establish an improved version of Theorem~E 
in the setting of Theorem~A?
\eprob

\bprob\label{HLP-prob2}
Can we establish the multidimensional versions of Theorems~A(c), B or C?
\eprob



The paper is organized as follows.  In Section \ref{HLP-sec2}, we present the main results of this paper.
In Theorem \ref{HLP-th1}, we present an affirmative answer to Problem \ref{HLP-prob1}.  
In Section \ref{HLP-sec4}, we state and prove five related theorems which extend three recent results of
Ponnusamy et al. \cite{KP201701,KP2018,PVW2019} from the case of analytic functions to the case of holomorphic functions
of several complex variables. In particular, Theorems \ref{HLP-th2} - \ref{HLP-th5} provide an affirmative answer to
Problem \ref{HLP-prob2}.

\section{Main Results}\label{HLP-sec2}

We first state the multidimensional version of Theorem~A(b). 

\bthm\label{HLP-th1}
If the series \eqref{Liu13} converges in the domain $Q$,  $p\in (0, 2]$, $a=|f(0)|<1$, and the estimate $|f(z)|<1$ holds in $Q$, then
\begin{equation}
|f(0)|^p+\sum_{k=1}^{\infty}|P_k(z)|+
\left(\frac{1}{1+a}+\frac{r}{1-r}\right)\sum_{n=1}^{\infty}|P_k(z)|^2\leq 1
\label{Liu15}
\end{equation}
for $z$ in the homothetic domain $r_p(a)Q$ and $r\leq r_p(a)$, where
$$r_p(a)=\frac{1-a^p}{2-a^2-a^p}.
$$
Moreover, if $Q$ is convex, then the numbers $r_p(a)$ and the factor $\frac{1}{1+a}$ in \eqref{Liu15} cannot be improved.
\ethm

\br
Note that $r_1(a)=1/(2+|f(0)|)$, $r_2(a)=1/2$ and $\ds \inf _{|a_0|<1} r_p(|a_0|) =p/(2+p)$ for  $0< p\leq 2$.
In particular, Theorem \ref{HLP-th1} for the cases $p=1$ and $p=2$ provide a multidimensional version of Theorem~A 
in an improved formulation of Theorem~E 
and multidimensional version of Theorem~A(b), 
respectively.
\er

%
%
%

Next, we state a generalization of Theorem~A(c). 

\bthm\label{HLP-th2}
Suppose that $Q$ is a complete circular domain centered at $0\in Q\subset\mathbb{C}^n$. If the series \eqref{Liu13}
converges in $Q$ such that $f(0)=0$ and  $|f(z)|<1$  for all $z\in Q$, then
\begin{equation}
\sum_{k=1}^{\infty}|P_k(z)|\leq 1
\label{Liu021}
\end{equation}
in the homothetic domain $\frac{1}{\sqrt{2}}Q$. Moreover, if $Q$ is convex, then the number $\frac{1}{\sqrt{2}}$ cannot be improved. 
\ethm

Note that the unit polydisk $\mathbb{D}^n$ is a convex complete circular domain centered at $0\in Q\subset\mathbb{C}^n$
and thus, Theorem \ref{HLP-th2} gives the following.

\bcor\label{HLP-cor3}
If the series \eqref{Liu17}  converges in the polydisk $\mathbb{D}^n$ such that $f(0)=0$ and  $|f(z)|<1$ for all $z\in \mathbb{D}^n$  then
\begin{equation}
\sum_{k=1}^{\infty}\bigg |\sum\limits_{|\alpha|=k}c_{\alpha}\, z^{\alpha}\bigg |\leq 1
\label{Liu022}
\end{equation}
in the polydisk $\frac{1}{\sqrt{2}}\mathbb{D}^n$. Moreover, the number $\frac{1}{\sqrt{2}}$ is sharp. 
\ecor

Let $K_n^0$ denote the largest number such that if the holomorphic function $f(z)=\sum\limits_{\alpha} c_{\alpha}\, z^{\alpha}$ satisfies $f(0)=0$, and $|f(z)|<1$ in the unit polydisk $\mathbb{D}^n$, then
\begin{equation}
\sum_{k=1}^{\infty}\sum\limits_{|\alpha|=k}|c_{\alpha}\, z^{\alpha}|\leq 1
\label{Liu023}
\end{equation}
holds in $K_n^0\cdot\mathbb{D}^n$. 

\bthm\label{HLP-th2-1}
It is true for $n> 1$ that 
$K_n^0\geq1/{\sqrt{2n}}$
and
$$
K_n^0\leq\left\{
\begin{array}{lll}
\ds \frac{1}{\sqrt{2}}&  ~\mbox{ for $n=2$},\\[4mm]
\ds \frac{2\sqrt{\log n}}{\sqrt{n}}& ~\mbox{ for $n>2$}.
\end{array}
\right.
$$
\ethm

Now we state a couple of generalizations of Theorem~B. 

\bthm\label{HLP-th3}
Assume that the series \eqref{Liu13} converges in the domain $Q$ and  $|f(z)|<1$ for all $z\in Q$. Then for $p\in (0,2]$, we have
\begin{equation}
|f(z)|^p+\sum_{k=N}^{\infty}|P_k(z)|
\leq 1
\label{Liu23}
\end{equation}
in the homothetic domain $(R_{N,p})Q$, where $R_{N,p}$ is the positive root of the equation
\begin{equation}
2(1+r)r^{N}-p(1-r)^{2}=0.
\label{Liu24}
\end{equation}
Moreover, if $Q$ is convex, then the radius $R_{N,p}$ is best possible. In particular,
\be\label{eq1-extra}
R_p:=R_{1,p}=\frac{p}{\sqrt{4p+1}+p+1}.
\ee
\ethm


Finally, we state the multidimensional version of  Theorem~C 
in a more general setting.

\bthm\label{HLP-th5}
Suppose that the series \eqref{Liu13} converges in the domain $Q$ and  $|f(z)|<1$ for all $z\in Q$. Also, let  $a=|f(0)|<1$ and
$p>0$.  Then
\begin{equation}
|f(z)|^p+\sum_{k=1}^{\infty}|P_k(z)|+ \left(\frac{1}{1+a}+\frac{r}{1-r}\right) \sum_{n=1}^{\infty}|P_k(z)|^2\leq 1
\label{Liu27a}
\end{equation}
for $z$ in the homothetic domain $r_{a,p}\cdot Q$ and $r\leq r_{a,p}$, 
where $r_{a,p}$ is the minimum positive root in $(0, 1)$ of the equation \eqref{Liu34}.  
Moreover, if $Q$ is convex, then the number $r_{a,p}$ cannot be improved.
\ethm

In the case of $0<p\leq 2$ of Theorem \ref{HLP-th5}, one can obtain the radius independent of the initial coefficient $a_0$ in a precise form. Because
of its independent interest, we may state it here explicitly.

\bcor
Suppose that the series \eqref{Liu13} converges in the domain $Q$ and  $|f(z)|<1$ for all $z\in Q$. Also, let  $a=|f(0)|<1$ and
$p\in (0,2]$.  Then \eqref{Liu27a} holds in the homothetic domain $R_p\cdot Q$, where $R_p$ is given by \eqref{eq1-extra}.
Moreover, if $Q$ is convex, then the number $R_{p}$ cannot be improved.
\ecor

The cases $p=1$ and $p=2$ in Theorem \ref{HLP-th5} give  multidimensional versions of \eqref{Liu015}  in Theorem~C, 
respectively. More precisely, we derive the following.

\br
\begin{enumerate}
\item For $p=1$, it is a simple exercise to see from Theorem \ref{HLP-th5} that
$$r_{a,1}= \frac{2}{3+a+\sqrt{5}(1+a)} ~\mbox{ and }~ \inf _{a\in [0,1)}r_{a,1} = \sqrt{5} -2.
$$
\item For $p=2$, we observe from Theorem \ref{HLP-th5} that $r_{a,2}=r_a$, where $r_a$ is the unique positive root of the equation
$$(1-a^3)r^3-(1+2a)r^{2}-2r+1=0.
$$
Also, we see that $1/3<r_a<1/(2+a)$, and $\inf _{a\in [0,1)}r_{a,2}=1/3$.
\end{enumerate}
\er
%
%
%
%


\section{Key lemmas and their Proofs}\label{HLP-sec3}
In order to establish our main results, we need the following lemmas. 
Our first lemma is a generalization of Theorem~B. 

\blem\label{lem2:Theo-C} \, 
Suppose that $f\in \mathcal{B}$ and $f(z)=\sum_{n=0}^{\infty}a_{n} z^{n}$. Then  for $p\in (0,2]$ and $N\in \IN$, we have the sharp inequality:
\begin{equation}
|f(z)|^p+B_N(f,r) \leq 1~ \mbox{ for }~ r\leq R_{N,p},
\label{Liu15a}
\end{equation}
where $R_{N,p}$ is the positive root of the equation \eqref{Liu24}.
The radius $R_{N,p}$ is the best possible.
\elem
\bpf
For functions $f \in\mathcal{B}$, the classical Schwarz-Pick lemma gives that
$$|f(z)|\leq \frac{r+a}{1+r a}, \quad |z|\leq r,\quad a=|f(0)| \in [0,1).
$$
Moreover, as $f \in\mathcal{B}$, we have $|a_n|\leq 1-|a|^2$ for all $n\geq 1$ and thus it follows that
\be\label{LP-eq2}
|f(z)|^p+B_N(f,r)  \leq  \left (\frac{r+a}{1+r a}\right )^p+ (1-a^2)\frac{r^N}{1-r}  = 1-\Psi_{N,p}(a)
\ee
where
\be\label{Liu35b}
\Psi_{N,p}(a)= \frac{1-r-(1-a^2)r^N}{1-r} -\left (\frac{r+a}{1+r a}\right )^p, \quad a \in [0,1].
\ee
Now, we wish to determine conditions such that $\Psi_{N,p}(a)\geq 0$ for all $a\in [0,1]$.
Note that $\Psi_{N,p}(1)=0$. We claim that $\Psi_{N,p}$ is a decreasing function of $a$, under the conditions of the
theorem. A direct computation shows that
$$\Psi_{N,p}'(a)=\frac{2ar^N}{1-r} -p(1-r^2)\frac{(r+a)^{p-1}}{(1+r a)^{p+1}}
$$
and
$$\Psi_{N,p}''(a)=\frac{2r^N}{1-r} -p(1-r^2)\frac{(r+a)^{p-2}}{(1+r a)^{p+2}} [p-1-2ar-(p+1)r^2].
$$
Evidently, $\Psi_{N,p}''(a)\geq 0$ for all $a\in [0,1]$, whenever $0<p\leq 1$. Hence for $r\leq R_{N, p}$,
$$\Psi_{N,p}'(a) \leq 
\Psi_{N,p}'(1)=\frac{2r^N}{1-r} -p\left (\frac{1-r}{1+r}\right )\leq 0,
$$
by the assumption that \eqref{Liu24} holds. Thus for each $r\leq R_{N, p}$ and $0<p\leq 1$, $\Psi_{N,p}$ is a decreasing function of $a\in [0,1]$, which implies that
$\Psi_{N,p}(a)\geq \Psi_{N,p}(1)=0$ for all $a\in [0,1]$ and the desired inequality \eqref{Liu15a} follows from \eqref{LP-eq2}.

Next, we show that condition $\Psi_{N,p}'(1)\leq 0$ is also sufficient for the function  $\Psi_{N,p}(a)$ to be decreasing on $[0,1]$ in the
case when $1<p\leq 2$. To do this we introduce an auxiliary function
$$ \Phi(r)= (1+r)^2\frac{(r+a)^{p-1}}{(1+r a)^{p+1}}, ~~r\in [0,1),
$$
and show that $\Phi$ is an increasing function of $r$ in $[0,1)$. An elementary computation leads to
$$ \Phi'(r)= (1+r)\frac{(r+a)^{p-2}}{(1+r a)^{p+2}}A(r), ~~r\in [0,1),
$$
where
\begin{eqnarray*}
A(r)&=&2(1+r a)(r+a) +(1+r)(p-1)(1+r a)-a(p+1)(1+r)(r+a)\\
&=& (1-a)\big [r\big (1-a+p(1+a)\big ) +a(p+1)+p-1\big ].
\end{eqnarray*}
Thus, for $p>1$ and $a\in [0,1]$, it follows that $A(r)\geq 0$ for all $r\in [0,1)$, and therefore we deduce that
$\Phi'(r)\geq 0$ for all $r\in [0,1)$. Hence,
$$\Phi(r)\geq \Phi(0)=a^{p-1} ~\mbox{ for all $r\in [0,1)$ and for $a\in [0,1)$.}
$$
This observation helps to derive that for $r\leq R_{N, p}$,
\begin{eqnarray*}
\Psi_{N,p}'(a)&= &\frac{2ar^N}{1-r} -p\left (\frac{1-r}{1+r}\right ) \Phi(r)\\
&\leq &
a^{p-1}\left [\frac{2a^{2-p}r^N}{1-r}-p\left (\frac{1-r}{1+r}\right ) \right ]\\
&\leq& a^{p-1}\left [\frac{2r^N}{1-r}-p\left (\frac{1-r}{1+r}\right ) \right ]=  a^{p-1} \Psi_{N,p}'(1)\leq 0 ,
\end{eqnarray*}
since $0\leq a^{2-p}\leq 1$ for $1<p\leq 2$. Again, $\Psi_{N,p}(a)$ is a decreasing function of $a\in [0,1]$, whenever $1<p\leq 2$
which implies that $\Psi_{N,p}(a)\geq \Psi_{N,p}(1)=0$ for all $a\in [0,1]$ and thus, the desired inequality \eqref{Liu15a} holds.

To prove that the radius is sharp, we consider the function $\varphi_a \in\mathcal{B}$ given by
$$\varphi_a(z)=\frac{a-z}{1-a z}=a-(1-a^2)\sum\limits_{k=1}^\infty a^{k-1}z^k,\quad a\in [0,1).
$$
For this function, with $a_0=a$ and $a_k=(a^2-1)a^{k-1}$ for  $k\in\mathbb{N}$, direct computation yields
\beq\label{KayPon8-eq1}
|\varphi_a(-r)|^p+B_N(\varphi_a,r) &=&\left (\frac{r+a}{1+r a}\right )^p+ (1-a^2)\frac{a^{N-1}r^N}{1-a r}\nonumber \\
&=& 1+\frac{(1-a)Q_{N,p}(a,r)}{(1+ar)^p(1-ar)}
\eeq
where
$$
Q_{N,p}(a,r)= (1-ar)(1+ar)^p\left [\left (\frac{1+a}{1-ar}\right ) a^{N-1}r^N -\frac{1}{1-a}\left (1-   \left (\frac{r+a}{1+r a}\right )^p\right )\right ],
$$
and it is easy to see that the last expression on the right of \eqref{KayPon8-eq1} is bigger than or equal to $1$ if and only if $Q_{N,p}(a,r)\geq 0$. In fact,
for $r> R_{N,p}$ and $a$ close to $1$, we see that
$$
\lim_{a\rightarrow 1^{-}} Q_{N,p}(a,r)=(1-r)(1+r)^p \left [\frac{2r^N}{1-r} -p\left (\frac{1-r}{1+r}\right )\right ]>0,
$$
showing  that the number $R_{N,p}$ in \eqref{Liu15a} is best possible. The proof of Lemma \ref{HLP-lem4} is done.
\epf

\br
Lemma \ref{lem2:Theo-C} for $p=1,2$ was obtained in \cite{KP201701}, see Theorem~B. 
Moreover, $R_{1,1}= \sqrt{5}-2$ and  $R_{1,2}=1/3$.
\er

In the case of $N=1$ of Lemma \ref{lem2:Theo-C}, the following refined formulation holds in sharp form.

\blem\label{lem3:Theo-C} \, 
Suppose that $f\in \mathcal{B}$ and $f(z)=\sum_{n=0}^{\infty}a_{n} z^{n}$ with $a=|f(0)|$ and $f_0(z)=f(z)-f(0)$. Then for $p\in (0,2]$, we have the sharp inequality:
\begin{equation*}
|f(z)|^p+B_1(f,r)+\left(\frac{1}{1+a}+\frac{r}{1-r}\right)\|f_0\|_r^2 \leq 1~ \mbox{ for }~ r\leq R_{p}=\frac{p}{\sqrt{4p+1}+p+1}.
\end{equation*}
The radius $R_{p}$ is the best possible.
\elem
\bpf
From the classical Schwarz-Pick lemma and Theorem~A(a), 
it follows that
$$
|f(z)|^p+B_1(f,r)+\frac{1}{1+a}\left (\frac{1+ar}{1-r}\right )\|f_0\|_r^2 \leq  \left (\frac{r+a}{1+r a}\right )^p+ (1-a^2)\frac{r}{1-r}  = 1-\Psi_{p}(a),
$$
where $\Psi_{p}(a)=\Psi_{1,p}(a)$, and $\Psi_{1,p}(a)$ is obtained from \eqref{Liu35b} by setting $N=1$ so that
$$\Psi_{p}(a)= \frac{1-(2-a^2)r}{1-r} -\left (\frac{r+a}{1+r a}\right )^p, \quad a \in [0,1].
$$
The rest of the proof follows from the proof of Lemma \ref{lem2:Theo-C}.
\epf

For our purpose, it is natural to ask about Lemma \ref{lem3:Theo-C} for $p>2$. Finding the radius $R_p$, independent of the constant term
$a_0=f(0)$, seems to be tedious. However, it is possible to state it in the following form, which is a generalization of Theorem~C. 

\blem\label{HLP-lem4}
Suppose that $f\in \mathcal{B}$ has the expansion $f(z)=\sum_{n=0}^{\infty} a_{n} z^{n}$, $a=|a_0|$, and $p>0$. Then
\begin{equation}
|f(z)|^p+B_1(f,r)+\left(\frac{1}{1+a}+\frac{r}{1-r}\right)\|f_0\|_r^2\leq 1
\label{Liu33}
\end{equation}
for $|z|=r\leq r_{a,p}$, where $r_{a,p}$ is the minimum positive root in $(0, 1)$ of the equation $A_{a,p}(r)=0$. Here $A_{a,p}(r)$ is given by \eqref{Liu34}
and the radius $r_{a,p}$ is best possible.
\elem
\bpf
Proceeding exactly as in the proof of Lemma \ref{lem2:Theo-C}, we find that
\begin{eqnarray}
|f(z)|^p+B_1(f,r)+\frac{1}{1+a}\left (\frac{1+ar}{1-r}\right )\|f_0\|_{r}^2
 &\leq & \left (\frac{r+a}{1+r a}\right )^p+\frac{r}{1-r}(1-a^2)\nonumber\\
 &=&1-\frac{A_{a,p}(r)}{(1-r)(1+a r)^p},
\label{Liu35}
\end{eqnarray}
where $A_{a,p}(r)$ is given by \eqref{Liu34}. Note that $A_{a,p}(r)$ is related by
$$\Psi_{1,p}(a)=\frac{A_{a,p}(r)}{(1-r)(1+a r)^p}
$$
where $\Psi_{1,p}(a)$ is obtained from \eqref{Liu35b} by setting $N=1$.
Evidently (\ref{Liu33}) holds if $A_{a,p}(r)\geq 0$.

Since $A_{a,p}(0)=1-a^p>0$ and $A_{a,p}(1)=(a^2-1)(1+a)^p<0$, we obtain that $A_{a,p}(r)\geq 0$ holds if and only if $r\leq r_{a,p}$, where $r_{a,p}$ is
the minimum positive root in $(0, 1)$ of the equation (\ref{Liu34}), i.e., $A_{a,p}(r)=0$. This gives that (\ref{Liu33}) holds for $r\leq r_{a,p}$.

To prove that the radius is sharp, we consider the function $\varphi \in\mathcal{B}$ given by
$$
\varphi(z)=\frac{a-z}{1-a z}=a-(1-a^2)\sum\limits_{k=1}^\infty a^{k-1}z^k,\quad a\in [0,1).
$$
and let $\varphi_0(z)=\varphi (z)-\varphi (0)$.  For this function, with $a_0=a$ and $a_k=(a^2-1)a^{k-1}$ for $k\in\mathbb{N}$, direct computations yields
$$ |\varphi(-r)|^p+B_1(\varphi,r)+\frac{1}{1+a}\left (\frac{1+ar}{1-r}\right )\|\varphi_0\|_r^2
= \Big(\frac{r+a}{1+r a}\Big)^p+ (1-a^2)\frac{r}{1-r}.
$$
Comparison of this expression with the right hand side of the expression in the formula (\ref{Liu35}) delivers the asserted sharpness.
The proof of Lemma \ref{HLP-lem4} is complete.
\epf

\section{Proofs of the main results}\label{HLP-sec4}

First we recall the following lemma due to Djakov and Ramanujan \cite[Lemma 1]{DjaRaman-2000}.

\begin{Lem}\label{HLP-lem3} (\cite{DjaRaman-2000})
For each holomorphic function $f:\,\mathbb{D}^n\to \mathbb{D}$ with $f(z)=\sum\limits_{\alpha} c_{\alpha}z^{\alpha}$, we have
\begin{enumerate}
\item[{\rm (a)}] $\ds b_k:=\bigg (\sum\limits_{|\alpha|=k}|c_{\alpha}|^2 \bigg)^{1/2}\leq 1-|c_0|^2$ for $k=1,2,\ldots$;

\item[{\rm (b)}] $\ds \sum\limits_{k=1}^{\infty}(b_{k})^{q}\leq (1-|c_0|^2)^{q-1}$ for $q\geq 2$.
\end{enumerate}
\end{Lem}

\subsection{Proof of Theorem \ref{HLP-th1}}
In each section of the domain $Q$ by the complex line
\begin{eqnarray}
\Lambda =\{z=(z_1,\ldots, z_n):\,z_j=a_j t,\, j=1,\ldots, n,\, t\in\mathbb{C}\},
\label{liu41}
\end{eqnarray}
the series turns into the power series in the complex variable $t$:
\begin{eqnarray}
f(at)=\sum_{k=0}^{\infty}P_k(a)t^k=f(0)+\sum_{k=1}^{\infty}P_k(a)t^k.
\label{liu42}
\end{eqnarray}
Since $|f(at)|<1$ for all $t\in\mathbb{D}$, by Theorem~A(b), 
we have
\begin{eqnarray}
|f(0)|^p+\sum_{k=1}^{\infty}|P_k(a)t^k|+\left (\frac{1}{1+|f(0)|}+\frac{r}{1-r}\right )\sum_{n=1}^{\infty}|P_k(a)t^k|^2\leq 1
\label{liu43}
\end{eqnarray}
for $z$ in the section $\Lambda\cap \left (\frac{1-a^p}{2-a^2-a^p}\cdot Q\right )$ and $r\leq \frac{1-a^p}{2-a^2-a^p}$, where $a=|f(0)|$.
The last inequality is just \eqref{Liu15}, since $\Lambda$ is an arbitrary complex line passing through the origin.

Moreover, if $Q$ is convex, then $Q$ is an intersection of half-spaces
\begin{eqnarray*}
Q=\bigcap\limits_{a\in J}\{z=(z_1,\ldots, z_n):\,{\rm Re}(a_1z_1+\cdots +a_n z_n)<1\}
\label{liu44}
\end{eqnarray*}
with some $J$. Because $Q$ is circular, we have
\begin{eqnarray}
Q=\bigcap\limits_{a\in J}\{z=(z_1,\ldots, z_n):\, |a_1z_1+\cdots +a_n z_n|<1\}.
\label{liu45}
\end{eqnarray}

Now it is sufficient to show that the constant $r_p(a)=\frac{1-a^p}{2-a^2-a^p}$ cannot be improved for each domain
$P_a=\{z=(z_1,\ldots, z_n):\, |a_1z_1+\cdots +a_n z_n|<1\}$.

In fact, for $a\in [0,1)$, there exists a function
$$\varphi_a(z_1)=\frac{a-z_1}{1-a z_1}=\sum_{n=0}^{\infty}a_n z_1^n\,\, (z_1\in\mathbb{D}, \, a_0=a,\, a_n=(a^2-1)a^{n-1},\, n=1,2,\ldots)
$$
such that $|\varphi_a(z_1)|<1$ in the unit disk $\mathbb{D}$, but for any $|z_1|=r>r_p(a)$,
$$a^p+B_1(\varphi,r)+  \left(\frac{1}{1+a}+\frac{r}{1-r}\right)\|\varphi_0\|_r^2 =a^p+ (1-a^2)\frac{r}{1-r}>1,
$$
so that  the inequality in (b) in Theorem~A 
fails in the disk $\mathbb{D}_r$. Thus we may finish the proof by using the function $f(z)=\varphi_a(a_1z_1+\cdots +a_n z_n)$. \hfill $\Box$
\vskip 2mm

\subsection{Proof of Theorem \ref{HLP-th2}}

By means of Theorem~A(c) 
and the analogous proof of Theorem \ref{HLP-th1}, we may easily verify that $\sum_{k=1}^{\infty}|P_k(z)|\leq 1$
in the homothetic domain $(1/\sqrt{2})Q$.

Now we prove that if $Q$ is convex, then the number $1/\sqrt{2}$ cannot be improved.
In fact, using the analogous proof of Theorem \ref{HLP-th1}, it is sufficient to show that the constant $1/\sqrt{2}$ cannot be improved for each domain $P_a=\{z=(z_1,\ldots, z_n):\,|a_1z_1+\cdots +a_n z_n|<1\}$.

Indeed, for $a=1/\sqrt{2}$, there exists a function
$$\psi(z_1)=z_1\frac{a-z_1}{1-a z_1}=\sum_{n=1}^{\infty}a_n z_1^n\quad \,\, (z_1\in\mathbb{D}, \, a_1=a,\, a_n=(a^2-1)a^{n-2},\, n=2,3,\ldots)
$$
such that $|\psi(z_1)|<1$ in the unit disk $\mathbb{D}$ and $\psi(0)=0$, but for any $|z_1|=r>1/\sqrt{2}$,
$$
\sum_{n=1}^{\infty}\left|a_{n}\right| r^{n}=\left . a r+ (1-a^2)\frac{r^2}{1-ar}\right |_{a=1/\sqrt{2}} = \frac{r/\sqrt{2}}{1-(r/\sqrt{2})}>1
$$
which implies that conclusion in Theorem~A(c) 
fails in the disk $\mathbb{D}_r$. Thus we may finish the proof  of Theorem \ref{HLP-th2} by using the function $f(z)=\psi(a_1z_1+\cdots +a_n z_n)$. \hfill $\Box$ 

\subsection{Proof of Theorem \ref{HLP-th2-1}}


If $z\in r\cdot \mathbb{D}^n$, for each holomorphic function $f(z)=\sum\limits_{\alpha} c_{\alpha}\, z^{\alpha}$ satisfying $f(0)=0$ and $|f(z)|<1$ in $\mathbb{D}^n$,  
it follows from Lemma~F 
and $c_0=f(0)=0$ that
\begin{eqnarray*}
\sum\limits_{k=1}^{\infty}\sum\limits_{|\alpha|=k}|c_{\alpha}|^2\leq 1-|c_0|^2=1.
\end{eqnarray*}

Therefore, using Cauchy-Schwarz inequality and above inequality,
we obtain that 
\begin{eqnarray*}
\sum_{k=1}^{\infty}\sum\limits_{|\alpha|=k}|c_{\alpha}\, z^{\alpha}|&\leq &\sqrt{\sum\limits_{k=1}^{\infty}\sum\limits_{|\alpha|=k}|c_{\alpha}|^2}\sqrt{\sum\limits_{k=1}^{\infty}\sum\limits_{|\alpha|=k}|z^{\alpha}|^2}
\leq  \sqrt{\sum\limits_{k=1}^{\infty} \left (\sum\limits_{j=1}^n |z_j|^2\right )^k}\\
&\leq & \sqrt{\sum\limits_{k=1}^{\infty} (n r^2)^k}=\sqrt{\frac{n r^2}{1-n r^2}}.
\end{eqnarray*}
We see that (\ref{Liu023}) holds if the last quantity is less than or equal to $1$. This gives the condition $r\leq \frac{1}{\sqrt{2n}}$.
Hence $K_n^0\geq \frac{1}{\sqrt{2n}}$.

Now we prove $K_n^0\leq \frac{1}{\sqrt{2}}$ for $n=2$. In fact, we consider the function
$$
f_0(z)=z_1\frac{a-z_2}{1-a z_2}=a z_1-(1-a^2)\sum_{k=1}^{\infty}a^{k-1}z_1z_2^k,
$$
where $a=1/\sqrt{2}$. Then for $z=(z_1,z_2)\not\in \frac{1}{\sqrt{2}}\cdot\overline{\mathbb{D}^2}$ with $|z_1|+|z_2|>\sqrt{2}$, we have
\begin{eqnarray*}
\sum_{k=1}^{\infty}\sum\limits_{|\alpha|=k}|c_{\alpha}\, z^{\alpha}|&= & a |z_1| + \sum_{k=2}^{\infty}(1-a^2)a^{k-2}|z_1||z_2|^k\\
&= & a |z_1| + (1-a^2)|z_1|\frac{|z_2|}{1-a|z_2|}=\frac{|z_1|}{\sqrt{2}-|z_2|}> 1.
\end{eqnarray*}
This implies $K_2^0\leq \frac{1}{\sqrt{2}}$. Thus the proof of the theorem follows from Theorem~D. 
\hfill $\Box$

\subsection{Proof of Theorem \ref{HLP-th3}}

By means of (\ref{Liu15a}) of Lemma \ref{lem2:Theo-C}, using the analogous proof of Theorem \ref{HLP-th1}, we may verify that
\begin{equation*}
|f(z)|^p+\sum_{k=N}^{\infty}|P_k(z)|\leq 1
\end{equation*}
in the homothetic domain $(R_{N,p})Q$, where $R_{N,p}$ is the positive root of the equation (\ref{Liu24}).

Now we prove that if $Q$ is convex, then the number $R_{N,p}$ cannot be improved.

In fact, using the analogous proof of Theorem \ref{HLP-th1}, it is sufficient to show that the constant $R_{N,p}$ cannot be improved for
each domain $P_a=\{z=(z_1,\ldots, z_n):\, |a_1z_1+\cdots +a_n z_n|<1\}$.

Indeed, from Lemma \ref{lem2:Theo-C}, it follows that for each $r>R_{N,p}$, there exists a function $f_1(z)=\sum_{n=0}^{\infty}a_{n} z^{n}$ such that $|f_1(z)|<1$ in $\mathbb{D}$, but (\ref{Liu15-1}) fails in the disk $\mathbb{D}_r$. Thus we may finish the proof of Theorem \ref{HLP-th3} by using the function $f(z)=f_1(a_1z_1+\cdots +a_n z_n)$.

Finally, for $N=1$, \eqref{Liu24} reduces to  $(2-p)r^2 +2(1+p)r-p=0$. Solving this gives the root $R_p$ given by\eqref{eq1-extra}
and the proof is complete.
\hfill $\Box$

\subsection{Proof of Theorem \ref{HLP-th5}}

By means of Lemma \ref{HLP-lem4}, using the analogous proof of Theorem \ref{HLP-th1}, we may verify that
\begin{equation*}
|f(z)|^p+\sum_{k=1}^{\infty}|P_k(z)|+\left (\frac{1}{1+a}+\frac{r}{1-r}\right )\sum_{n=1}^{\infty}|P_k(z)|^2\leq 1
\end{equation*}
for $z$ in the homothetic domain $r_{a,p}\cdot Q$ and $r\leq r_{a,p}$, where $r_{a,p}$ is the minimum positive root in $(0, 1)$ of the equation (\ref{Liu34}). Now we prove that if $Q$ is convex, then the number $r_{a,p}$ cannot be improved. The reasoning as in the proof of Theorem \ref{HLP-th3} concludes the proof of  Theorem \ref{HLP-th5}. \hfill $\Box$

\subsection*{Acknowledgments}
The work of the first author is supported by Guangdong Natural Science Foundation (Grant No. 2021A030313326)
and the work of the second author is supported by Mathematical Research Impact Centric Support (MATRICS) of
the Department of Science and Technology (DST), India  (MTR/2017/000367).
The authors thank the referees very much for their valuable comments and suggestions to this paper.

%

\end{document}